\documentclass[a4paper,10pt]{amsart}

\usepackage[all]{xy}
\usepackage{amssymb}

\newcommand{\Spe}{{\rm Spec}}

\newcommand{\A}{{ \rm Aut }}

\newcommand{\Z}{{\mathbb{Z}}}

\renewcommand{\mod}{{\;\rm mod}}
\newtheorem{pro}{Proposition}[section]
\newtheorem{lemma}[pro]{Lemma}

\newtheorem{cor}[pro]{Corollary}
\newtheorem{theorem}[pro]{Theorem}
\newtheorem{definition}[pro]{Definition}

\title{On cyclic covers of the projective line}
\author{Antoniadis A. Jannis \and Kontogeorgis Aristides}

\email{antoniad@math.uoc.gr}
\address{
Department of Mathematics, University of Crete, 71409 Herakleion, Crete, Greece,
Greece\\ { \texttt{\upshape http://www.math.uoc.gr/\~ \;\!\!antoniad}}
}

\email{kontogar@aegean.gr}
\address{
Department of Mathematics, University of the \AE gean, 83200 Karlovassi, Samos,
Greece\\ { \texttt{\upshape http://eloris.samos.aegean.gr}}
}
\date{\today}

\begin{document}
\bibliographystyle{amsplain}

\maketitle

\begin{abstract}
We construct configuration spaces for cyclic covers of the projective line that admit extra automorphisms and we describe the 
locus of curves with given  automorphism group. As an 
application we provide examples of arbitrary high genus 
that are defined over their field of moduli and are not hyperelliptic.
\end{abstract}
\section{Introduction}
In \cite{Gut-Sha} J. Gutierrez and T. Shaska introduced the notion of Dihedral invariants 
in order to study moduli of hyperelliptic curves with extra involutions.
A hyperelliptic curve $X$  is a Galois cover of order $2$ of the projective line and corresponds 
to a Galois extension of order $2$ of the rational function field. Let $j$ be the generator 
of the Galois group  $Gal(X/\mathbb{P}^1_k)$, {\em i.e.} the hyperelliptic involution.
It is known that $\langle j \rangle$ is a normal subgroup of the whole automorphism 
group $\A (X)$, so the {\em reduced group} $\bar{G}=\A(X)/\langle j \rangle$ makes sense 
and if the genus $g$ of $X$  is $g\geq 2$, then the reduced group 
  is a finite subgroup of $\A (\mathbb{P}^1_k)=PGL_2(k)$. 
In this setting a curve $X$ has an extra involution if there  is an involution in the 
automorphism group that does not induce the identity in the reduced group.

Let $\mathcal{L}_g$ denote the locus of hyperelliptic curves with extra involutions.
Gutierrez and Shaska were able to 
classify hyperelliptic curves that admit extra involutions. Moreover the space 
$\mathcal{L}_g$ was proved to be rational \cite[th. 3.4]{Gut-Sha}. Finally they applied their theory in 
order to study the classical problem  of field of moduli and of field of definition.

For $n \in \mathbb{N}$ let $C_n$ denote the cyclic group of order $n$.
As a generalization of the above, we consider $n$-cyclic covers $X$  of the projective line of genus $g_X \geq 2$.
In this case it is no more true that the cyclic group $C_n=Gal(X/\mathbb{P}^1_k)$ is a normal 
subgroup of the whole automorphism group. Under some mild hypotheses 
we can assume the normality and then we consider 
in this case also the reduced automorphism group 
$\bar{G}:=\A (X)/\langle \sigma \rangle$, where $\sigma$ is a generator of the cyclic group 
$Gal(X/\mathbb{P}^1_k)$. 
Under this normality assumption 
the group of automorphisms  was studied by  R. Brandt, H. Stichtenoth 
 \cite{Bra-Sti},\cite{Brandt-PhD} if $n$ is a prime number 
and by the second author  if $n$ is
 composite \cite{Ko:99}. 

 Let $F$ be the the function field of the curve $X$  and 
let $F_0$ be the rational function field, $F_0=F^{C_n}$.
In the general case by an {\em extra automorphism} we mean a non 
trivial  automorphism of order $\delta$, $\delta | n$, $\delta\neq 1$ that is not an 
element of the cyclic group $\mathrm{Gal}(F/F_0)$.   
Let $\mathcal{L}_{n,s,\delta}$ denote the locus of branched $n$-cyclic covers of $\mathbb{P}^1(k)$ branched 
over $s$-points  that admit an extra automorphism of order $\delta \mid n$. 

We are working over a field $k$ of characteristic $p\geq 0$. If $p>0$ then we also assume that $(p,n)=1$. 
We define two $n$-covers of the projective line that admit an extra automorphism  to be equivalent if they 
are isomorphic and if the isomorphism transforms the extra 
automorphism  of the first curve 
to a generator of the cyclic group generated by the extra automorphism of the second curve. Using this stronger notion of equivalence we 
are able to 
 construct  configuration spaces for the set of Galois cyclic  covers 
$X\rightarrow \mathbb{P}^1$ with $Gal(X/\mathbb{P}^1_k)=C_n$ such that 
$C_n \lhd \A (X)$. This is a generalization of the theory of hyperelliptic curves where $n=2$. 
The usage of a  stronger notion of equivalence is a common practice in the study of 
moduli spaces of curves \cite{HarrisModuli}, where $n$-pointed curves or level structures on the 
Jacobians are introduced.

Moreover we are able to describe the locus of curves with reduced automorphism 
group isomorphic to a dihedral group and we can provide examples 
of curves with reduced group isomorphic to a given finite automorphism group of 
the projective line. In all cases the structure of the automorphism group is given.

In the next section we give all possible automorphism groups for cyclic covers 
of the projective line that have genus $\leq 10$. Our methods are applicable to higher 
genera as well, but the length of the results avoid us for giving them here.

Finally, using these configuration spaces we are able to 
provide examples of curves of arbitrary high genus, 
 that are defined over their field of moduli and are not 
hyperelliptic.  The problem of providing rational models over the 
field of moduli is a difficult task, and  as far as the authors know 
for $g>3$ the only examples known are for hyperelliptic curves \cite{Sha:hyper2004},\cite{Gut-Sha}.
\section{Computations}
Let $k$ be an algebraic closed field of characteristic $p \geq 0$.
Let $F_0=k(x)$ be the function field of the projective line $\mathbb{P}^1(k)$. We consider 
a cyclic extension of $F_0$ of the form $F:=k(x,y)$ where 
\begin{equation} \label{defeqfun}
y^n=\prod_{i=1}^s (x-\rho_i)^{d_i}=:f(x), \;\; 0<d_i<n.
\end{equation}

If $d:=\sum_{i=1}^s d_i \equiv 0 \mod n $ then the place at infinity does not ramify at the 
above extension \cite[p.667]{Ko:99}. The only places of $F_0$ that are ramified are 
the places $P_i$ that correspond to the points $x=\rho_i$ and the corresponding 
ramification indices are given by  
\[
e_i=\frac{n}{(n,d_i)}.
\] 
Moreover if $(n,d_i)=1$ then  the places $P_i$ are ramified completely and the 
Riemann-Hurwitz formula implies that the function field  $F$ has genus
\begin{equation} \label{eq2}
g=\frac{(n-1)(s-2)}{2}.
\end{equation}
Notice that the condition $g\geq 2$ is equivalent to $s \geq 2 \frac{n+1}{n-1}$. In particular, $s \geq 2$.
\begin{pro}
Let $G=\A(F)$. 
Suppose that a cyclic extension $F/F_0$ of the rational function field  $F_0$ is 
ramified completely at $s$ places and $n:=|\mathrm{Gal}(F/F_0)|$. If $2n<s$ then 
$\mathrm{Gal}(F/F_0) \lhd G$.
\end{pro}
\begin{proof}
\cite[prop. 1]{Ko:99}
\end{proof}

From now on we will assume that  $\mathrm{Gal}(F/F_0) \lhd G$.

\begin{lemma} \label{2.2}
Let $\gamma \in G$ and $\bar{\gamma}$ be its image in $\bar{G}$. Suppose that 
$\bar{\gamma}$ is an element of order $\delta |n$. The elements $\gamma,\bar{\gamma}$ have the 
same order if and only if $\bar{\gamma}$ does not fix any of the places $\rho_i$. 
\end{lemma}
\begin{proof}
\cite[Th 9]{Ko:99}
\end{proof}

Fix a generator $\sigma$ of $\mathrm{Gal}(F/F_0)$. 
We notice that if $\tau$ is an extra automorphism then the set 
$\{\tau \sigma^i$, $i=0,\ldots,n\}$ consists of extra automorphisms.  

\begin{lemma} \label{lemma-norm}
Suppose that $\tau$ is an extra automorphism of $F$, and
let $s$ be the number of ramified places 
at the  extension $F/F_0$ and let $d$ be the degree of the defining polynomial.
 Then $\delta \mid s$, $\delta\mid d$ and  the defining equation of $F$ can be written 
as 
\begin{equation} \label{eq*}
y^n= \sum_{i=0}^{d/\delta} a_i  x^{\delta \cdot i}, 
\end{equation}
where $a_0=1$.
\end{lemma}
\begin{proof}
Notice that $s=d$ if and only if $d_i=1$ for all ramified places.
Since, by lemma \ref{2.2},  $\bar{\tau}$ has the same order with $\tau$, $\bar{\tau}$  permutes the $s$-points $P_i$ without 
fixing them, thus $\delta:=\mathrm{ord}(\tau) \mid s$. On the other hand, 
we can choose the coordinate $x$ so that $\bar{\tau}(x)=\tau(x)=\zeta x$, where 
$\zeta$ is a primitive $\delta$-root of one. We can change the enumeration so that
the ramified points are of the form
$\{ \zeta^i \beta_j,$ $i=0,\ldots,\delta-1,j=1,\ldots,s/\delta\}$. Therefore, the 
defining equation of the curve is of 
the form (we allow multiple roots $\beta_i$)
\[y^n=\prod_{i=1}^{d/\delta}( x^\delta -\beta_i^\delta)=\sum_{i=0}^{d/\delta} a_i  x^{\delta \cdot i}.
\]
 The coefficients  $a_0,\ldots,a_{d/\delta}$ 
are  the symmetric polynomials of $\beta_i^\delta$, multiplied by $(-1)^{\frac{s}{\delta}-i}$. 

In particular the constant term is  
\begin{equation} \label{zero-coeff}
a_0=(-1)^{d/\delta} \prod_{i=1}^{d/\delta} \beta_i ^\delta. 
\end{equation}
and by changing   $x$ to $\lambda x$ for suitable $\lambda$ we can 
 assume that the constant term $a_0=1$.  
\end{proof}
\begin{definition}
We will say that the $n$-cover of the projective line is in  {\em normal 
form} if and only it is given by an equation:
\begin{equation} \label{normal-form}
y^n=x^{s} + \sum_{i=1}^{\frac{d}{\delta}-1} a_i x^{\delta i} + 1.
\end{equation}
\end{definition}

In what follows we will make the assumption that $s=d$, {\em i.e.} all polynomials defining normal forms 
have simple roots. We define the spaces
$
\mathcal{N}_{n,s,\delta}
$  to be the spaces of normal forms corresponding to $n$-th cyclic covers of the projective line 
ramified above $s$ points and with an extra automorphism of order $\delta$. 
Each space $\mathcal{N}_{n,s,\delta}$ is a Zariski open set in  $\mathbb{A}_k$ corresponding to
the complement $V(D)$, where $D$ is a suitable discriminant and itself an algebraic variety 
with coordinate ring the ring $k[x_1,\ldots,x_{s-2}]_D$.  
The assumption on the 
simplicity of the roots of the polynomials corresponding to the normal forms, implies that 
all curves corresponding to points in $\mathcal{N}_{n,s,\delta}$ have the same genus. 

\section{Dihedral Invariants}
A natural candidate for parametrizing the set of $n$-covers of the 
projective line that admit an extra automorphism of order $\delta$, is the set of coefficients 
$\{a_{s/\delta-1},\ldots, a_1\}$ of a normal form. But two different normal forms 
might correspond to isomorphic curves, therefore  
we have to
see in what extend this normalization determines the variable $x$. 
The condition $\tau(x)=\zeta x$, implies that $\bar{\tau}$ fixes the  places  $0,\infty$. Moreover 
we can change the defining equation  by a morphism  $\gamma \in \mathrm{PGL}(2,k)$ of the form 
$\gamma:x \rightarrow m x$ or $\gamma:x \rightarrow \frac{m}{x}$ so that the new equation is 
again in normal form.
By applying $\gamma$ to  (\ref{zero-coeff}) we obtain that 
\[
(-1)^{s/\delta} \prod_{i=1}^{s/\delta}  \gamma(\beta_i)^\delta=1,
\] 
and this gives  $m^s=(-1)^{s/\delta}$. 
As in the case of hyperelliptic curves, studied by 
Gutierrez and Shaska \cite{Gut-Sha},  we obtain that $x$ is determined up to 
a coordinate change by the subgroup $D_{s/\delta}$ generated by 
\[
\tau_1:x \rightarrow \varepsilon x, \;\;\;\tau_2:x \rightarrow \frac{1}{x},
\]
where $\varepsilon$ is a primitive $s/\delta$-root of one. 

The action of $D_{s/\delta}$ on the parameter space $k(a_1,\ldots,a_{s/\delta})$ is given by 
\[
\tau_1:a_i \mapsto \varepsilon ^{\delta i} a_i, \mbox{ for } i=1,\ldots,s/\delta 
\]
\[
\tau_2:a_i \mapsto a_{d/\delta-i}, \mbox{ for } i=1,\ldots, [ s/ 2\delta ]. 
\]
Notice that if $s/\delta=1$ then the above actions are trivial, therefore 
the normal form determines the equivalence class. 
If $s/\delta=2$ then 
\[
\tau_1(a_1)= -a_1,\;\; \tau_1(a_2)=a_2 \;\; \tau_2=\mathrm{Id},
\]
and the action is not dihedral but cyclic on the first factor.

In what follows, in order to have a dihedral action  we will assume that $s/\delta > 2$. 
We have the following:
\begin{pro}
Assume that $s/\delta > 2$.
The fixed field $k(a_1,\ldots,a_{s/\delta})^{D_{s/\delta}}$ is the same as the function field of the 
variety $\mathcal{L}_{n,s,\delta}$.
\end{pro}

\begin{lemma}
Let $r:=s/\delta>2$.
The elements
\begin{equation} \label{defui}
u_i:=a_1^{r-i}a_i+a_{r-1}^{r-i}a_{r-i}, \;\;\mbox{ for }  1 \leq i \leq r, 
\end{equation}
are invariant under the action of the group $D_{s/\delta}$ defined above.
The elements $u_i$ are called the dihedral invariants of $D_{s/\delta}$.
\end{lemma}
\begin{theorem}
Let $u=(u_1,\ldots,u_{r})$ be the $r$-tuple of dihedral invariants. Then $k(\mathcal{L}_{s,n,\delta})=k(u_1,\ldots,u_r)$. 
\end{theorem}
\begin{proof}
The dihedral invariants are fixed by the $D_r$ action. Hence $k(u)\subseteq k(\mathcal{L}_{s,n,\delta})$ therefore 
we have  $ [k(a_1,\ldots,a_r):k(u)] \geq 2r$.
We will show that $[k(a_1,\ldots,a_r):k(u)]=2r=|D_r|$. For each $2 \leq i \leq r-2$ we 
have 
\[
a_1^{r-i}a_i +a_{r-1}^{r-i}a_{r-i}=u_i
\]
\[
a_1^i a_{r-i} + a_{r-1}^i a_i =u_{r-i}.
\]
The above two equations imply that $a_i,a_{r-i} \in k(u,a_1,a_{r-1})$.
Moreover
\begin{equation} \label{1122}
u_1=a_1^r+a_{r-1}^{r}\;\;\mbox{ and } u_{r-1}=2a_1 a_{r-1},
\end{equation}
and, by substitution, we arrive at
\[
2^r a_{r-1}^{2r} - 2^r a_{r-1}^r u_1 + u_{r-1}^r=0,
\]
so $a_{r-1}$ satisfies a  polynomial of degree $2r$ 
over $k(u)$. On the other hand (\ref{1122}) implies that 
$a_1 \in k(u,a_{r-1})$ so  $[k(a_1,\ldots,a_r):k(u)]\leq 2r$ therefore  $[k(a_1,\ldots,a_r):k(u)]=2r$.
\end{proof}

\begin{definition}
Let $(X_a,\tau)$, $(X_{a'},\tau')$ be two $n$-covers of the projective 
live, ramified completely over $s$-points and let $\tau$  (resp. $\tau'$) be  an 
extra automorphism  of $X_a$ (resp. $X_{a'}$). We say that these two 
covers are isomorphic if there is an isomorphism between the corresponding function fields $F_a,F_a'$,
$\alpha:F_a \rightarrow F_{a'}$  such that $\alpha \langle \tau \rangle  \alpha^{-1} =\langle \tau' \rangle$.
\end{definition}

\begin{lemma}
Let $1<\delta$ be a divisor of $n$ and consider  
$a:=(a_1,\ldots,a_r) \in k^r$. Equation (\ref{normal-form}) defines a $n$-branched cover $X_a$  of the projective line 
such that the reduced automorphism group contains an extra automorphism $\tau$ of the form 
$x \mapsto  \zeta x$, where $\zeta$ is a primitive $\delta$-root of one. 
   Two such pairs $(X_a,\tau)$, $(X_{a'},\tau')$ are isomorphic if and only if the corresponding 
   dihedral invariants are the same. 
\end{lemma}
\begin{proof}
It is trivial to see that equations of the form (\ref{normal-form}) indeed define curves with an  extra automorphism 
$\tau$  so that $\tau(x)=\zeta x$ and conversely every $n$-branched cover of  the projective 
line after a suitable choice of the model can be written in such a way. 

Suppose that two pairs $(X_a,\tau)$, $(X_{a'},\tau')$ are isomorphic. Let us 
denote by $F_a,F_{a'}$ the corresponding function fields.  By definition 
there is an isomorphism $\alpha: F_a \rightarrow F_{a'}$ such that $\alpha\tau \alpha^{-1}=(\tau')^\nu$, $(\nu,\delta)=1$. 
 We observe first 
that $F_a$ (reps. $F_{a'}$) has a unique rational subfield $k(x)$ (resp. $k(x')$) of degree $n$.
Therefore $\alpha(x) \in k(x')$. The conjugation condition $\alpha \tau \alpha^{-1}=(\tau')^\nu$, implies 
that $\alpha(x) = {\frac{\lambda}{\lambda_3}} x'$. Indeed, since $\alpha$ is an invertible morphism from $k(x)\rightarrow k(x')$ it
is of the form $a(x)=\frac{\lambda x' + \lambda_1}{\lambda_2 x' +\lambda_3}$ with 
\[
\begin{pmatrix} \lambda & \lambda_1 \\ \lambda_2 & \lambda_3 \end{pmatrix} \in PGL(2,k)
\]
and the commuting relation implies, $\lambda_1=\lambda_2=0$.
 The function fields  $F_a,F_{a'}$ are defined by  equations of the form:
\begin{equation} \label{Xa}
y^n = x^s+\sum_{i=1}^{r-1} a_i x^{i \delta}+1
\end{equation}
\begin{equation} \label{Xa1}
{y'}^{n}={x'}^s+ \sum_{i=1}^{r-1} a_i' {x'}^{i \delta}+1.
\end{equation}
Since the above two curves are isomorphic the coefficients $\{a_1,\ldots,a_r\}$ can be transformed 
to the coefficients $\{a_1',\ldots,a_r'\}$  under the group  generated by  $\tau_1, \tau_2$. This proves that the 
 curves $X_a,X_{a'}$  have the same dihedral invariants. 

Conversely if the curves defined by equations (\ref{Xa}) and (\ref{Xa1}) share the same 
dihedral invariants then they are isomorphic, and the isomorphism is given by $\alpha(x)=x'$, 
and $\alpha \tau \alpha^{-1} = \tau'$. 
\end{proof}

Let $\bar{G}$ be the reduced automorphism of the curve $X$ and let $\tau$ be an extra automorphism
of order $\delta,1<\delta, \delta \mid n$. The elements $g \tau g^{-1}$ are also extra automorphisms, and 
if $g \tau g^{-1} \neq \tau$ then $(X,\tau)$, $(X,g \tau g^{-1})$ correspond to different 
points in the space $\mathcal{L}_{d,n,\delta}$.  From the above remark we obtain:
\begin{pro}
The number of points in $\mathcal{L}_{n,s,\delta}$ that correspond to a  curve $X$ with an 
 an extra automorphism $\tau$ of order $\delta$, is equal to the number of conjugacy classes of 
the cyclic group generated by  $\tau$ in $G$. 
\end{pro}

\begin{cor}
If the reduced automorphism is Abelian, then to the curve $X$ with an extra automorphism
of order $\delta,1<\delta, \delta \mid n$ corresponds only one point in $\mathcal{L}_{n,s,\delta}$.
\end{cor}

\section{Equation of curves with extra automorphisms}
Consider the variety $\mathcal{L}_{n,s,\delta}$ parametrized by the set of dihedral invariants $u_1,\ldots,u_r$. 
We will try to characterize the locus of extra automorphisms. 
\begin{pro} \label{listcases}
If the characteristic $p$ of the field $k$ is zero or $p>5$ then 
the possible automorphisms groups of the projective line are given 
by the following list
\begin{enumerate}
\item Cyclic group of order $\delta$.
\item $D_\delta=\langle \sigma,\tau\rangle$, $(\delta,p)=1$ where 
$\sigma(x)=\xi x$, $\tau(x)=1/x$, $\xi$ is a primitive $\delta$-th 
root of one. The possible orbits of the $D_\delta$ action are 
$B_\infty=\{0,\infty\}$,
$B^-=\{\mbox{roots of } x^\delta-1\}$,
$B^+=\{\mbox{roots of } x^\delta+1\}$, 
$B_a=\{\mbox{roots of } x^{2\delta} + a x^\delta +1\}$, $a\in 
k\backslash \{\pm 2\}$.
\item 
$A_4=\langle \sigma,\mu\rangle$, $\sigma(x)=-x$, $\mu(x)=i\frac{x+1}{x-1}$, $i^2=-1$. The possible orbits of action are 
the following sets:
$B_0=\{0,\infty,\pm 1, \pm i\}$,
 $B_1=\{\mbox{roots of } x^4-2i \sqrt{3}x^2+1\}$, 
$B_2=\{\mbox{roots of } x^4+2i \sqrt{3}x^2+1\}$,
$B_a=\{\mbox{roots of } \prod_{i=1}^3 (x^4-a_ix^2+1)\}$, 
where $a_1\in k \backslash \{\pm 2, \pm 2i \sqrt{3}\}$, 
$a_2=\frac{2a_1+12}{2-a_1}$, $a_3=\frac{2a_1-12}{2+a_1}$
\item 
$S_4=\langle \sigma,\mu \rangle$, $\sigma(x)=ix$, $\mu(x)=i\frac{x+1}{x-1}$, $i^2=-1$. The possible orbits of action are 
the following sets:
$B_0=\{0,\infty,\pm 1, \pm i\}$, $B_1=\{\mbox{roots of } 
x^8+14x^4+1\}$, 
$B_2=\{\mbox{roots of } (x^4+1)(x^8-34x^4+1) \}$, 
$B_a=\{\mbox{roots of } (x^8+14x^4+1)^3-a(x^5-x)^4\}$, $a\in 
k\backslash \{108\}$.
\item 
$A_5=\langle \sigma,\rho\rangle$, $\sigma(x)=\xi x$, 
$\rho(x)=-\frac{x+b}{bx+1}$, where $\xi$ is a primitive fifth root
of one and $b=-i(\xi+\xi^4)$, $i^2=-1$.
The possible orbits of action are 
the following sets:
$B_\infty=\{0,\infty\}\cup \{\mbox{roots of } f_0(x):=x\big(x^{10}+11 ix^5+1 \big)\}$, 
$B_0=\{\mbox{roots of } f_1(x):=x^{20}-228i x^{15}-494 x^{10} -228 ix ^5+1\}$,
$B_0^*=\{\mbox{roots of }
x^{30}+522 i x^{25}+ 10005(x^{20}-x^{10})-522ix^5-1\}$,
$B_a=\{\mbox{roots of }  f_1(x)^3 - a f_0(x)^5\}$, 
where $a\in k\backslash\{ -1728i\}$
\item Semidirect products of elementary Abelian groups with 
cyclic groups: 
$\big(\mathbb{Z}/p\mathbb{Z}\times \cdots \times
 \mathbb{Z}/p\mathbb{Z} \big) \rtimes \mathbb{Z}/m \mathbb{Z}$,
of order $p^tm$, where $m \mid p^t-1$. The possible orbits of action are
$B_\infty=\{\infty\}$, 
$B_0=\{\mbox{roots of } f(x):=x \prod_{j=1}^{\frac{p^m-1}{m}}(x^m-b_j)\}$
where $b_j$ are selected so that the elements of the  additive group 
$\mathbb{Z}/p\mathbb{Z}\times \cdots \times
 \mathbb{Z}/p\mathbb{Z}$ seen as elements in $k$, 
are roots of $f(x)$.
 $B_a=\{\mbox{roots of } f^m-a\}$, where $a \in k\backslash B_0$.
\item $PSL(2,p^t)=\langle \sigma,\tau,\phi\rangle$,
where $\sigma(x)=\xi^2x$, $\tau(x)=-\frac{1}{x}$, $\phi(x)=x+1$ and 
$\xi$ is a primitive $m=p^t-1$ root of one. The orbits of the action are
$B_{\infty}=\{\infty,\mbox{roots of } x^m-x\}$. 
$B_0=\{\mbox{roots of } (x^m-x)^{m-1}+1\}$,
$B_a=\{\mbox{roots of } 
\big(
(x^m-x)^{m-1} +1
\big) ^{\frac{m+1}{2}}-a(x^m-x)^{\frac{ m(m-1)}{2}}\}$
where $a\in k^*$.
\item $PGL(2,p^t)=\langle \sigma,\tau,\phi\rangle$,
where $\sigma(x)=\xi x$, $\tau(x)=\frac{1}{x}$, $\phi(x)=x+1$ and 
$\xi$ is a primitive $m=p^t-1$ root of one. The orbits of the action are
$B_{\infty}=\{\infty,\mbox{roots of } x^m-x\}$. 
$B_0=\{\mbox{roots of } (x^m-x)^{m-1}+1\}$,
$B_a=\{\mbox{roots of } 
\big(
(x^m-x)^{m-1} +1
\big) ^{m+1}-a(x^m-x)^{m(m-1)}\}$
where $a\in k^*$.
\end{enumerate}
\end{pro}
\begin{proof}
    Brandt PhD, Stichtenoth \cite{StiEssen}.
\end{proof}
{\bf Remark} Assume that the reduced group $\bar{G}$ is one of the above groups. For the action of $\bar{G}$ on the 
projective line we have the generic orbits $B_a$ where $\# B_a=|\bar{G}|$ and the 
special orbits where $\#B_a = |\bar{G}|/|\bar{G}(P)|<|\bar{G}|$ and $\bar{G}(P)$ is the stabilizer of a point in the 
special orbit $|B_a|$.  If the extra automorphism $\tau$ of order $\delta$ reduces to an element $\bar{\tau}$  in $G(P)$ for some $P$ in 
an orbit $B_a$ then there are $|\bar{G}|/|\bar{G}(P)|$-conjugates of $\bar{\tau}$ in $G$.

\subsection{Cyclic Case}
%-----------------------------------------
\begin{pro} \label{cyc-inside}
Let $F$ be a function field  corresponding to  a point in $\mathcal{L}_{d,n,\delta}$ 
and denote by $\tau$ the extra automorphism of order $\delta$. Assume that 
this automorphism is an element in the group $\langle \tau' \rangle$, where $\tau'$ is an extra automorphism of order 
$1<\delta'$ with $\delta \mid \delta'$, $\delta< \delta'$. Then the corresponding values of the dihedral invariants are $u_1=\cdots=u_r=0$. 
\end{pro}
\begin{proof}
Since there is an automorphism $\tau$ of order $\delta$ we can write a defining  equation for the curve $F$ 
of the form 
\[
y^n=f(x),
\]
 where the polynomial $f(x)$ can be written as 
\[
f(x)=\prod_{i=1}^{s/\delta} (x^\delta - \beta_i^\delta) = x^{s} + a_{r-1} x^{(\frac{d}{\delta}-1)\delta} + \cdots + a_1 x^\delta+a_0.
\]
The existence of an automorphism  $\tau'$ implies a further condition on the set of $\beta_i$'s  {\em i.e.} 
they should be of the form $\{  \varepsilon^k \beta_j, k=0,\ldots,\frac{\delta'}{\delta} \}$. This gives us that the polynomial 
$f(x)$ is of the form 
\[
f(x)=\prod_{j=1}^{s/\delta'} (x^{\delta'}-\beta_j^{\delta'}),
\] 
{\em i.e.}, $a_i=0$ if $\frac{\delta'}{\delta} \nmid i$. 
This gives that $a_1=a_{r-1}=0$ thus 
by (\ref{defui}) we have that all $u_i=0$.
\end{proof}
In order to assume that $u\neq 0$ we have to assume that the extra automorphism $\tau$ of order $\delta$ is not 
contained in any cyclic subgroup of the reduced automorphism group of greater order. 

Let $e=\frac{\delta'}{\delta}$.
We can change the dihedral invariants to the following invariants:
\[
u_i^{(e)}= a_e ^{r-1} a_i + a_{r-e} ^{r-1} a_{r-i}.
\]
These are also dihedral invariants  and, by using the same method, we can 
prove that they generate again  the rational function field $k(\mathcal{L}_{n,d,\delta})$. Moreover, in the new generators
the locus of curves with an extra automorphism of order $\delta'$ is 
given by $u^{(e)}_i=0$ for $e \nmid i$.

{\bf Remark 1.} The above situation is a usual phenomenon for birational 
maps. A birational map can collapse a subvariety to a point. We will come back to 
this on  remark 3 on page \pageref{dihedral-blowup}.

\begin{lemma} \label{normalformreduce}
Suppose that we have the $n$-cover given by the equation 
\[
C:y^n=\sum_{i=0}^r a_i x^{\delta i}, \mbox{ with }  a_{r}=1.
\]
Then, the coordinate change $(x,y)\mapsto (\lambda x,  a_0 ^{1/n} y)$ where $\lambda^{\delta r}=a_0$
makes the curve $C$ isomorphic to the curve 
\[
C':y^n=\sum_{i=0} ^r a_i' x^{\delta i}, \mbox{ with } a_{r}'=a_0'=1.
\]
\end{lemma}
\begin{proof}
We observe that after change $x\mapsto \lambda x$ the polynomial 
$f(x)=\sum_{i=0}^r a_i x^{\delta i}$ is transformed to the polynomial 
\[
f'(x)=\sum_{i=0}^r a_i \lambda^{\delta i } x^{\delta i}= 
\]
\[
=a_0\left( x^{r\delta}+\lambda^{(r-1)\delta} \frac{a_{r-1}}{a_0} x^{(r-1)\delta}+
 \cdots + \lambda^\delta\frac{a_1}{a_0} x^\delta +1 \right)
\]
and $a_0$ is absorbed by the corresponding element from $y$.

Let $a_i'$ denote the coefficient corresponding to $x^i$. By computation 
\begin{equation} \label{nfr2}
a_i'=\lambda^{\delta i} \frac{a_i}{a_0}=\lambda^{\delta(i-r)}{a_i}.
\end{equation}

\end{proof}
\subsection{Dihedral Case} 
%-----------------------------------
%
%
Let $D_m$ be the dihedral group generated by  $\tau,\sigma$ where $\tau^m=1$, $\sigma^2=1$ and $\sigma \tau \sigma^{-1}=\tau^{-1}$. 
We can choose our model so that 
  $\tau(x)=\zeta x$, where $\zeta$ is a primitive 
$m$-root of one, then $\sigma$ is of the form $\sigma(x)=a/x$ and 
the fixed points of $\sigma$ are $\pm \sqrt{a}$.
Moreover after a change of variable $x\rightarrow \frac{x}{a^{1/2}}$ we can assume that the dihedral group 
is given by $\tau(x)=\zeta x$, and $\sigma(x)=\frac{1}{x}$. The fixed points of $\sigma$ are $\pm 1$.

We distinguish two cases:

{\bf Case A.} In this case we assume that the extra automorphism of order $\delta$ inside the reduced group $D_m$ is an 
element of the cyclic group generated by $\tau$.

{\bf Case B.} In this case the extra automorphism is of order $2$, $2\mid n$ and the extra automorphism inside the 
reduced group $D_m$ is an element of order two, not inside the cyclic group generated by $\tau$.

We observe that if the reduced group is   the  abelian group $D_2$ the two cases coinside.
We begin our study with case A.
\begin{lemma} \label{dihed-sym}
Let $X: y^n =f(x)$ be a curve that admit  the dihedral group $G=D_\delta$ as reduced group. If the variable $x$ is chosen 
so that the dihedral group is generated by the above $\tau,\sigma$, then for the coefficients $a_i$ of the polynomial 
$f(x)=\sum_{i=0}^r a_i x^{i\delta}$ we have  $a_r=1$ and for $1\leq i \leq r-1$ we have $a_i = a_0 a_{r-i}$.
\end{lemma}
\begin{proof}
Since the reduced automorphism group $G$ has a cyclic subgroup of order $\delta$ the polynomial 
$f(x)$ is of the form $f(x)=\sum_{i=0}^r a_i x^{i\delta}$. Furthermore the reduced group permutes the 
roots of $f(x)$. In particular, if  $\beta_i$ is a root then $\frac{1}{\beta_i}$ is also a root. 
Let $ \zeta ^j \beta_i$ be the roots of $f(x)$. 

 Then by Vieta's formulae we have:
\[
a_i = (-1)^{r-i} \sum_{\nu_1,\ldots,\nu_i} \beta_{\nu_1}^\delta \cdots \beta_{\nu_i}^\delta=
(-1)^{r-i}\prod_{\nu=1}^r \beta_\nu^\delta
\sum_{\mu_1,\ldots,\mu_{r-i}}  \frac{1}{\beta_{\mu_1}^\delta} \cdots \frac{1}{\beta_{\mu_{r-i}}^\delta}
\]
\[
=(-1)^{-i} a_0 \sum_{\mu_1,\ldots,\mu_{r-i}}  \frac{1}{\beta_{\mu_1}^\delta} \cdots \frac{1}{\beta_{\mu_{r-i}}^\delta}
= a_0  a_{r-i}.
\]
\end{proof}

\begin{pro} \label{dih-form}
If $D_\delta$ is a subgroup of the reduced automorphism  group then for the dihedral invariants we have 
$2^{r-2}u_1^2=u_{r-1}^r$.
\end{pro}
\begin{proof}
We observe first that if $\beta\neq 1$ is a root of the polynomial $f$ then $\frac{1}{\beta}$ is also a root of $f$ and their product 
is equal to one. We distinguish the following cases:

{\bf Case 1.} None of the fixed points of  $\sigma$  are included in the branch locus. 
Then, $r$ is even, $a_0=1$ and $a_1=a_{r-1}$. 
In this case $u_1=2a_1^r$ and $u_{r-1}=2a_1^2$ and we have that 
$u_{r-1}^{\frac{r}{2}}=2^{\frac{r}{2}}a_1^r=2^{\frac{r-2}{2}}u_1$, 
{\em i.e.,} the desired result.

{\bf Case 2.} From the fixed points of $\sigma$, only 
 $-1$ is included in the branch locus.
Then $r$ is odd, and $a_0=(-1)^{r+1}=1$ and the above
curve is in normal form so we proceed as above.

{\bf Case 3.} From the fixed points of $\sigma$ 
only  $1$ is included in the 
set of roots of $f$. In this case  $r$ is odd and 
$a_0=(-1)^r=-1$. We have to use lemma 
\ref{normalformreduce} in order to compute the 
coefficients $a_i'$ of a normal form. By equation (\ref{nfr2}) 
we have:
\[
a_1'=\lambda^{\delta(1-r)}a_1 
\mbox{ and } 
{a_1'}^r=\lambda^{r\delta(1-r)}a_1^r=(-1)^{1-r}a_1^r =a_1^r
\] 
and
\[
a_{r-1}'=\lambda^{\delta(r-1-r)}a_{r-1}
\mbox{ so  }
{a_{r-1}'}^r=\lambda^{-\delta r}a_{r-1}^r=(-1)a_{r-1}^r=
(-1)a_0^r a_1^r=a_1^r.
\]
Therefore $u_1=2a_1^r$ and 
\[
u_{r-1}=2a_1'a_{r-1}'=2 \lambda^{-\delta r}a_1a_{r-1}=2(-1)a_0a_1^2=
2a_1^2,
\]
and the desired result follows.

{\bf Case 4.}
Both fixed points of $\sigma$ $\{\pm 1\}$ are included
in the set of roots of $f$. In this case   $r$ is even and 
$a_0=(-1)^r (-1)=-1$ 
We compute
\[
a_1'=\lambda^{\delta(1-r)}a_1 
\mbox{ and } 
{a_1'}^r=\lambda^{r\delta(1-r)}a_1^r=(-1)^{1-r}a_1^r =-a_1^r
\] 
and
\[
a_{r-1}'=\lambda^{\delta(r-1-r)}a_{r-1}
\mbox{ so }
{a_{r-1}'}^r=\lambda^{-\delta r}a_{r-1}^r=(-1)a_{r-1}^r=
(-1)a_0^r a_1^r=-a_1^r,
\]
and the desired result follows.
\end{proof}
\begin{pro}
Under the assumption that there is an extra automorphism on the 
$C_n$-cover of $\mathbb{P}^1_k$, the reduced 
automorphism group is dihedral and the extra automorphism lies inside the 
the cyclic part of the dihedral reduced group, the possible structures for the 
full automorphism group are:
\[
G_1:=\langle  R,S | R^{2n}=1, S^\delta=1,RSR^{-1}=S^{-1}\rangle,
\]
\[
G_2:=\langle R,S | R^{2n}=1,S^\delta=1,(RS)^2=1 \rangle,
\]
\[
G_3:=\mathbb{Z}/n\mathbb{Z} \rtimes D_\delta.
\]
The group $G_1$ appears in case 4 of proposition \ref{dih-form}, the group 
$G_2$ appears in cases 3,4 of proposition \ref{dih-form} and the group 
$G_3$ appears in case 1 of proposition \ref{dih-form}.
\end{pro}
\begin{proof}
\cite[Th. 15]{Ko:99}
\end{proof}
{\bf Remark 2.} In the case of even $r$ the equation $2^{r-1}u_1^2=u_{r-1}^r$ is reducible. The component 
$u_{r-1}^{r/2}=2^{r-2}u_1$ corresponds to the group $G_3$ and 
the component  $u_{r-1}^{r/2}=-2^{r-2}u_1$
corresponds to the group $G_1$. 
%%%%
% Blow up remark
%%%%

{\bf Remark 3.} \label{dihedral-blowup}
Let $X$ be a point in $\mathcal{L}_{d,n,\delta}$. Assume 
that the generator $\tau:x\mapsto \zeta x$, $\zeta^\delta=1$, 
is an element in the cyclic group $\langle \tau' \rangle$, 
where $\tau'$ is an extra automorphism of order 
$1<\delta < \delta'$, $\delta \mid \delta'$. 

By proposition  \ref{cyc-inside} the corresponding values 
of the dihedral invariants are $u_1=\cdots=u_r=0$.
If $e=\frac{\delta'}{\delta}$, then we consider the 
invariants 
\[
u_i^{(e)}=a_e ^{r-1}a_i+a_{r-e}^{r-1}a_{r-i}.
\]
Since $a_i$ are the coefficients of a normal form 
we have $a_0=1$ and lemma \ref{dihed-sym} implies
$a_{r-i}=a_i$, so 
\[
u_i^{(e)}=2 a_e^{r-1}a_i.
\]
Thus the ratio 
\[
\frac{u^{(e)}_i}{u_i}=\left(\frac{a_e}{a_1}\right)^{r-1},
\]
 is independed of $i$. 
Therefore,
\[
u_i^{(e)} u_j-u_j^{(e)}u_i=0,
\]
and the new invariants $u_1^{(e)},\ldots, u_r^{(e)}$, 
correspond to a blow-up at $(0,\ldots,0)$.

We now turn our attention to case B. In this case we assume that the extra automorphism reduces to an element $\sigma$ of order $2$
inside the reduced group $D_m$, $m>2$, and that $\sigma$ is not inside the cyclic subgroup of order $m$ in $D_m$.
We now have to use a different model for the action of automorphism group on the rational function field.

By diagonalizing the  $2\times 2$ matrices defining the action in the model given in  proposition \ref{listcases} we can assume that 
\[
\sigma:x \mapsto -x \;\;\;\;  \tau^\ell:x\mapsto -\frac{(\zeta^\ell+1)x-\zeta^\ell+1}{(\zeta^\ell-1)x-\zeta^\ell-1},\] 
where $\zeta$ is a primitive $m$-th root of unity.
The fixed points of $\sigma$ are $\{0, \infty\}$ and the fixed points of $\tau$ are $\pm 1$.
The orbit of an element $\beta$ that is not fixed by $\tau$ under the action of $\tau$ is 
$\left\{-\frac{(\zeta^\ell+1)\beta-\zeta^\ell+1}{(\zeta^\ell-1)\beta-\zeta^\ell-1},\ell=1,\ldots,m \right\}$.
Therefore the orbits of the action of the dihedral group is  one of the following:
\[
B_0:=\left\{\frac{1-\zeta^\ell}{1+\zeta^\ell},\ell=1,\ldots,m \right\}, 
\;\;\;
B_\infty:=\left\{\frac{1+\zeta^\ell}{1-\zeta^\ell},\ell=1,\ldots,m \right\}
\]
\[
B_{1}:=\{\pm 1\},\;\;\;\; 
B_\beta:=\left\{\mp\frac{(\zeta^\ell+1)\beta-\zeta^\ell+1}{(\zeta^\ell-1)\beta-\zeta^\ell-1},\ell=1,\ldots,m \right\},
\]
for $\beta \not \in \{\pm 1,0,\infty\}$.

Since we assume that there is an extra automorphism of order $2$ the orbits $B_0,B_\infty$ can not be in 
the ramification locus. So we can have $t_1$ orbits of the form $B_{\beta_i}$ and $t_0$ orbits of the form   $B_1$ in the ramification locus,
$t_0=0,1$.  
\begin{pro}
If $n\equiv 0 \mod 2$, the reduced automorphism group is isomorphic to $D_m$ and we have an extra automorphism 
of order $2$ that is not contained in the cyclic subgroup of $D_m$. The locus of such curves in $\mathcal{L}_{n,s,2}$ is of dimension $t_1$.
The full automorphism group is isomorphic to  $C_{nm}\rtimes C_2$ if $t_0=1$ and to $C_n \rtimes D_m$ if $t_0=0$.
\end{pro}
\begin{proof}
The locus $\mathcal{L}_{n,s,2}$ is of dimension $t_1$ since we have $t_1$ free variables $\beta_i$. The structure of the 
full automorphism group follows by \cite[th. 15]{Ko:99}
\end{proof}

In what follows we will give a natural construction for combining two cyclic $C_n$-covers $X_1,X_2$
of the projective line with reduced automorphism group $\bar{G}$, to a new curve $X$ that is also a 
a $C_n$-curve of the projective line and has reduced automorphism group $\bar{G}$, and 
the branch locus of $X$ is the union of the branch loci of $X_1,X_2$.

Let $\mathcal{N}_{n,s,\delta} \cong \mathbb{A}^{s/\delta-1}_k$ be the set 
of normal forms corresponding to $n$-covers of the projective line branched
over $s$-points, and admit an extra automorphism of order $\delta \mid n$.
For every point
 $a=(a_1,\ldots,a_{s/\delta-1}) \in 
\mathcal{N}_{n,s,\delta}$, 
we will denote by $f_a(x)\in k[x]$ the polynomial
\[
f_a(x)=x^{ s} + \sum_{i=1}^{s/r-1} a_i x^{\delta i} + 1.
\]
There is a natural map 
\begin{equation} \label{natural-forms-dec}
\Phi:\big( 
\mathcal{N}_{n,s_1,\delta} \times \mathcal{N}_{n,s_2,\delta}
\big)  \backslash \Delta \rightarrow  \mathcal{N}_{n,s_1+s_2,\delta}
\end{equation}
\[
(a,b)=\big( (a_1,\ldots,a_{s_1/\delta-1}),
(b_1,\ldots,b_{s_2/\delta-1})  \big) \mapsto 
(\gamma_1,\ldots,\gamma_{(s_1+s_2)/\delta-1}),
\]
where $1\leq \gamma_i\leq (s_1+s_2)/\delta-1$ correspond to the coefficients of the
polynomial $f_a \cdot f_b$, and  $\Delta$ corresponds to  the closed algebraic set given by the
resultant 
$\mathrm{Res}(f_a,f_b)=0$.
\begin{lemma}
If we are working over an algebraically closed field $k$,
the map $\Phi$ is onto.
\end{lemma}
\begin{proof}
The map $\Phi$ is onto since every polynomial of degree $s_1+s_2$ 
can be written as the product of two polynomials of corresponding degrees
$s_1,s_2$.
\end{proof}
\begin{lemma}
Let $r_\mu=s_mu/\delta$ $\mu=1,2$, and consider the dihedral automorphisms 
$\tau_1(a_i)=\varepsilon^{i\delta} a_i$, $\tau_2(a_i)=a_{r-i}$.
The map $\Phi$ is compatible with the $\tau_1,\tau_2$ action.
\end{lemma}
\begin{proof}
Consider the spaces $\mathcal{N}_{n,a,\delta}$, $a=s_1,s_2,s_1+s_2$. 
The generators of the dihedral action $\tau_{1,a},\tau_{2,a}$ act on them as follows:
\[
\begin{array}{ccc}
\tau_{1,s_1}(a_i)=\varepsilon^{i\delta} a_i &\tau_{1,s_2}(b_i)=\varepsilon^{i\delta} b_i &\tau_{1,s_1+s_2}(\gamma_i)=\varepsilon^{i\delta}\gamma_i, \\
\tau_{2,s_1}(a_i)= a_{\frac{s_1}{\delta}-1} &\tau_{2,s_2}(b_i )= b_{\frac{s_2}{\delta}-i}&\tau_{1,s_1+s_2}(\gamma_i)=\gamma_{\frac{s_1+s_2}{\delta}-i}, 
\end{array}.
\]
We have to prove that 
\[
\tau_{\mu,s_1+s_2}\circ \Phi =\Phi \circ (\tau_{\mu,s_1},\tau_{\mu,s_2}), \;\;\;\mu=1,2.
\]
This is true, for instance:
\[
\sum_{\kappa+\lambda=i} \tau_{2,s_1} (a_\kappa) \tau_{2,s_2}(b_\lambda)=\sum_{\kappa+\lambda=i} a_{\frac{s_1}{\delta}-\kappa} \beta_{\frac{s_2}{\delta}-\lambda}=
\]
\[
=\sum_{\kappa+\lambda=\frac{s_1+s_2}{\delta}-i} a_\kappa b_\lambda= \tau_{2,s_1+s_2} (\gamma_i). 
\]
\end{proof}
$\mathbf{A_4}$ {\bf Case.} 
It is known and easy to check, that every element of  $A_4$ is conjugate to an element of order  $2$ or $3$.
As in the study of the dihedral group we have to distinguish two cases.

{\bf Case A}. In this case the reduced group is $A_4$ and the extra automorphism is of order two.
We will use the notation of proposition \ref{listcases}.
The components corresponding to the orbit $B_i$, $i=1,2$  are given by 
$\pm 2 i \sqrt{3} \in \mathcal{N}_{n,4,\delta}$.
Since the reduced automorphism $\sigma:x \rightarrow \zeta x$, where $\zeta$ is 
a primitive $\delta$-root of one, has to be 
lifted to an element of order $\delta$ in 
$G$ the orbit $B_0$ of $A_4$ could not 
be contained in the set of branch points 
of the cover $C \rightarrow \mathbb{P}^1_k$.
Set $t_1=1$ (resp. $t_2=1$) if the  orbit $B_1$ (resp. $B_2$) is contained in the set of branch points and $t_1=0$ 
(reps. $t_2=0$) if not. Let $t_3$ be the 
number of orbits of the form $B_a$ that 
are contained in the set of branch points.
Then  the number $s$ of branch points of the above cover is 
\[
s=4(t_1+t_2)+12t_3.
\]
Moreover the locus  in 
$\mathcal{N}_{n,s,\delta}$ is the image 
of the loci of $A_4$ in $\mathcal{N}_{n,12,\delta}$ and $\mathcal{N}_{n,4,\delta}$ under successive 
applications of the appropriate $\Phi$ functions.
 
By computation, 
the component in $\mathcal{N}_{n,12,\delta}$ corresponding to 
an orbit $B_a$ is given by 
\begin{eqnarray*}
a_1=& 
-{\frac {a \left( a-6 \right)  \left( a+6 \right) }{ \left( 2+a
 \right)  \left( -2+a \right) }}
\\
a_2=&
-33\\
a_3=&
2\,{\frac {a \left( a-6 \right)  \left( a+6 \right) }{ \left( 2+a
 \right)  \left( -2+a \right) }}\\
a_4=&
-33\\
a_5=&
-{\frac {a \left( a-6 \right)  \left( a+6 \right) }{ \left( 2+a
 \right)  \left( -2+a \right) }}
\end{eqnarray*}
Observe that $V_4=D_2$ is a normal subgroup of $A_4$ hence the symmetry
of the coefficients of $a_i$.
The dihedral invariants are now computed:
\[
u_1=2a_1^6,u_2=-66a_1^4,u_3=-4a_1^4,u^4=-66a_1^2,u_5=2a_1^2,
\]
and the locus of curves with $A_4$ as reduced 
group in 
$\mathcal{L}_{n,12,\delta}=\mathbb{A}^{5}$ 
is given by the 1-dimensional algebraic set 
with equations
\[
71874 u_1^2=u_2^3,u_3=\frac{2}{33}u_1, 
u_4^2=-66 u_2, -33u_5=u_4.
\]

It is very complicated to write down 
the locus in $\mathcal{L}_{n,s,\delta}$
of curves that admit $A_4$ as a reduced 
group, although for fixed values $t_1,t_2,t_3$ 
it can be done with the aid of a program of symbolic 
computation.

\begin{pro}
 The dimension of the locus of curves with reduced group   isomorphic to  
$A_4$ is $t_3$.  For the structure of the group $G$ of  automorphisms of the corresponding curves we have 
\[
G=
\left\{ 
\begin{array}{ll}
\mathbb{Z}/n\mathbb{Z} \times A_4  & \mbox{if } t_1=t_2=0  \\
\mathbb{Z}/3n\mathbb{Z} \times V_4 & \mbox{otherwise}
\end{array}
\right.
\]
\end{pro}
\begin{proof}
The dimension formula is clear since on the formula of the normal 
forms we have $t_3$ free variables. The structure theorem follows by theorem 
18 in \cite{Ko:99}. Notice that only cases a,b of theorem 
18 in \cite{Ko:99} can appear and that the 
action of $A_4$ in $\mathbb{Z}/n\mathbb{Z}$ is trivial.
\end{proof}

{\bf Case B.}
In this case the reduced automorphism group is again $A_4$ but the extra automorphism is 
of order $3$. Using a diagonalization argument we can change the coordinate $x$ for the rational function 
field, so that the group $A_4$ is generated by:
\[
\tau_1:x \mapsto 
\left(-\frac{1}{2}+\frac{1}{2}\sqrt{3}\right) x,
\;\;
\tau_2:x\mapsto 
 \frac{ x + 1+i}{(1-i)x-1}.
\]
The elements $\tau_1,\tau_2$ are of orders $3$ and $2$ respectively. The fixed points of $\tau_1$ are $\{0,\infty\}$
and the fixed points of $\tau_2$ are
\[\left\{ \left( 1/2+1/2\,i \right)  \left( 1+\sqrt {3} \right) , 
 \left( -1/2-1/2\,i \right)  \left( -1+\sqrt {3} \right)  \right\}.\]

We have the following orbits of $A_4=\langle \tau_1,\tau_2 \rangle$ acting on $\mathbb{P}^1(k)$.
\[
B_0=\left\{
0, 
1-\sqrt {3},
{\frac { \left( -1+i+\sqrt {3}+i\sqrt {3} \right)  \left( -1+\sqrt {3}
 \right) }{1-i+\sqrt {3}+i\sqrt {3}}},
-{\frac { \left( -1+i+\sqrt {3}+i\sqrt {3} \right) ^{2} \left( -1+
\sqrt {3} \right) }{ \left( 1-i+\sqrt {3}+i\sqrt {3} \right) ^{2}}}
\right\}
\] 
\[
B_\infty=\left\{\infty,
 {\frac {\sqrt {3}}{3+\sqrt {3}}},
-{\frac { \left( -1+i+\sqrt {3}+i\sqrt {3} \right) \sqrt {3}}{ \left( 
1-i+\sqrt {3}+i\sqrt {3} \right)  \left( 3+\sqrt {3} \right) }},
{\frac { \left( -1+i+\sqrt {3}+i\sqrt {3} \right) ^{2}\sqrt {3}}{
 \left( 1-i+\sqrt {3}+i\sqrt {3} \right) ^{2} \left( 3+\sqrt {3}
 \right) }}
\right\}
\]
\[
B_1:=\left\{
1,-1/2+1/2 i \sqrt{3},
i \sqrt{3}-1/2 \sqrt{3}+1-3/2 i,
\sqrt{3}-2,
\right.\]
\[\left.
-i \sqrt{3}-1/2 \sqrt{3}+1+3/2 i,-1/2-1/2 i \sqrt{3}
\right\}
\]
\begin{equation} \label{eq12}
B_\alpha:=\mbox{roots of the polynomial}
 f_a(x):=\left(
 {x}^{3}-{\alpha}^{3} \right)  \left( {x}^{3}-3\,{\frac {\sqrt 
{3} \left( \alpha-1+\sqrt {3} \right) ^{3}}{ \left( 3\,\alpha+\sqrt {3
}\alpha-\sqrt {3} \right) ^{3}}}
 \right)\cdot\end{equation}
\[\cdot  \left(
 {x}^{3}-3\,{\frac {
\sqrt {3} \left( -\alpha+i\alpha+\sqrt {3}\alpha+i\alpha\,\sqrt {3}-2-
4\,i+2\,i\sqrt {3} \right) ^{3}}{ \left( 6\,i\alpha+2\,\sqrt {3}\alpha
+4\,i\alpha\,\sqrt {3}+\sqrt {3}-i\sqrt {3}+3+3\,i \right) ^{3}}}
 \right) 
\cdot\]\[\cdot 
 \left( {x}^{3}-3\,{\frac {\sqrt {3} \left( i\alpha-\sqrt {3}
\alpha+i\sqrt {3}+3-i-\sqrt {3} \right) ^{3}}{ \left( -3\,\alpha-i
\sqrt {3}-3\,\sqrt {3}\alpha+3\,i\alpha+i\alpha\,\sqrt {3}-3 \right) ^
{3}}} \right) 
\]
where $\alpha$ is not an element in $B_0,B_\infty,B_1$.

Since there is an extra automorphism of order $3$ we have that the orbits $B_0,B_\infty$ 
could not be in the ramification locus.
\begin{theorem}
Assume that we have a $n$-th cover of the projective line that admits an extra automorphism 
of order $3$ so that the reduced group is of order $3$. Assume that we have $t_1$ orbits of the form $B_1$ and 
$t_2$ orbits of the form $B_\alpha$, $t_1=0,1$. The dimension of the locus of curves with reduced automorphism 
group  $A_4$ is equal to $t_2$, and for the structure of the full automorphism group $G$ we have:
\[
G=
\left\{
\begin{array}{ll}
\mathbb{Z}/n\mathbb{Z}\times A_4 & \mbox{ if } t_1=0\\
     \langle R,S | R^2=S^2,S^{2n}=1,RSR^{-1}=S^r \rangle                                & \mbox{ otherwise} 
\end{array}
\right.
\]
for a suitable $r$ that can can be explicitly described as solution of some modular equations.
\end{theorem}
\begin{proof}
The dimension argument is as before while the structure comes from 
\cite[th. 18]{Ko:99}.
\end{proof}

Before we proceed to the next case we will need the following 
\begin{lemma} \label{chang-coord}
Let $f(x)$ be a polynomial so that the roots of $f$ form orbits of a group $G \subset PGL_2(k)$, 
{\em i.e.} 
\[
f(x)=\prod_{g\in G} (x -g (\rho)), \quad \rho \in k.
\]
Let $\mu(x)=\frac{ax+b}{cx+d}$ be a M\"obius transformation.
Then an orbit of the element $\mu^{-1}\rho$ is given by the roots of the polynomial
\[
f\left(\mu(x) \right)\cdot (cx+d)^{\deg f}.
\] 
\end{lemma}
\begin{proof}
This is a direct computation with the rational function $f(\mu(x))$.
\end{proof}
%
%-----------------------------------------------------------------------
$\mathbf{S_4}$ {\bf Case.} 
We will use again the  the notation of proposition \ref{listcases}. All elements in 
$S_4$ are conjugate to one of the following elements  $\sigma,\mu,\sigma\mu^2$ with corresponding orders
$4,3,2$. In particular all elements in $S_4$ have order $2,4$ or $3$.

We distinguish the following cases:

{\bf Case A}
In this case 
we assume that the extra automorphism of order $\delta$ in the reduced automorphism group $S_4$ is in the cyclic group of order 
four generated by the automorphism $\sigma(x)=ix$. Then,
the orbit $B_0$ can 
not be included in the ramification locus, by the assumption that the 
curve has an extra automorphism of order $\delta$. Let $t_1=1$ 
(resp. $t_2$) 
$B_1$ (resp. $B_2$) is included in the ramification locus and $t_1=0$
(reps. $t_2=0$) if not. Let $t_3$ be the number of $B_a$ orbits.  
The number of branch points is 
\[
s=8t_1+6t_2+24t_3
\]
\begin{pro}
If we assume that the extra automorphism of order $\delta$ in the the reduced group $S_4$ is in the cyclic group 
generated by $\sigma(x)=ix$ then 
the locus of curves with  reduced group 
isomorphic to $S_4$ is of dimension $t_3$. The automorphism 
group in this case is $C_n \times S_4$.
\end{pro}
\begin{proof}
The dimension formula is clear, since we have $t_3$ free variables 
on the formula of the normal form. The structure of the group follows 
by \cite[Th. 20]{Ko:99}
\end{proof}
{\bf Example:} Using Maple we compute the dihedral invariants of the curve:
\[
y^n=\left(  \left( {x}^{8}+14\,{x}^{4}+1 \right) ^{3}-a \left( {x}^{5}-x
 \right) ^{4} \right)  \left(  \left( {x}^{8}+14\,{x}^{4}+1 \right) ^{
3}-b \left( {x}^{5}-x \right) ^{4} \right).
\]
\[
\begin{array}{ll}
u_1= &2\, \left( 84-b-a \right) ^{12} \\
 u_2= &      2\, \left( 84-b-a \right) ^{10} \left( 2946-38\,b-38\,a+ab \right) \\
u_3 =&2\, \left( 84-b-a \right) ^{9} \left( 55300-429\,b-429\,a-8\,ab
 \right) \\
u_4 =&2\, \left( 84-b-a \right) ^{8} \left( 588015-712\,b+28\,ab-712\,a
 \right) \\
u_5 =&2\, \left( 84-b-a \right) ^{7} \left( 3392424+7342\,b+7342\,a-56\,ab
 \right) \\
u_6= &2\, \left( 84-b-a \right) ^{6} \left( 8699676-12324\,b-12324\,a+70\,ab
 \right) \\
u_7= &2\, \left( 84-b-a \right) ^{5} \left( 3392424+7342\,b+7342\,a-56\,ab
 \right) \\
u_8= &2\, \left( 84-b-a \right) ^{4} \left( 588015-712\,b+28\,ab-712\,a
 \right) \\
u_9= &2\, \left( 84-b-a \right) ^{3} \left( 55300-429\,b-429\,a-8\,ab
 \right) \\
u_{10}= &2\, \left( 84-b-a \right) ^{2} \left( 2946-38\,b-38\,a+ab \right) \\
u_{11}= &2\, \left( 84-b-a \right) ^{2}\\
u_{12}=  &2\\
\end{array}
\]
{\bf Case B}
In this case we assume that the reduced group is still $S_4$ but the reduced automorphism of order $\delta$ is of order $3$.
Using a diagonalization argument we find that  for $Q(x)=-1/3\,{\frac { \left(  \left( 3+\sqrt {3} \right) x-\sqrt {3}+3
 \right) \sqrt {3}}{1-i- \left( 1-i \right) x}}$, the elements 
\[\sigma':=Q^{-1} \circ \sigma \circ Q=-{\frac { \left( -1+i+\sqrt {3}+i\sqrt {3} \right) x}{1-i+\sqrt {3}+i
\sqrt {3}}}
\]
\[Q^{-1} \circ \mu \circ Q= {\frac { \left( 3+3\,i+i\sqrt {3}-\sqrt {3} \right) x-i\sqrt {3}-3+3\,
i+\sqrt {3}}{ \left( -3+3\,i+i\sqrt {3}-\sqrt {3} \right) x-i\sqrt {3}
+3+3\,i+\sqrt {3}}}.
   \] 
On this new coordinates $\sigma'$ is an element of order $3$ fixing $0,\infty$.
The orbits of the action of $S_4$ in the new coordinates follow by applying $Q^{-1}$ on the orbits of the old coordinates.
In particular, since the lower left entry of the matrix $Q$ is zero  we compute   the following orbits:
\[
B_0':=\big\{ Q^{-1}(P), P\in \{0,\infty,\pm 1,\pm i\}
\big\}  
\]
\[
B_1':=\{\mbox{roots of } g_1(Q(x)), g_1(x):=x^8+14x^4+1\},
\]
\[
B_2':=\{\mbox{roots of } g_2(Q(x)), g_2(x):=(x^4+1)(x^8+34x^4+1) \}
\]
\[
B_a':=\{\mbox{roots of } f_a(Q(x)), f_a(x):=(x^8+14x^4+1)^3-a (x^5-x)^4 \}.
\]
Observe that $\{0,\infty\}$ are in $B_1'$ (by direct copmutation or since this orbit has $8=24/3$ elements, 
therefore this orbit can not be in the ramification locus.
\begin{pro}
Assume that the reduced automorphism group is $S_4$ and that the extra automorphism of order $\delta$ is of order $3$.
Then if we have $t_0,t_2,t_a$ orbits of type $B_0',B_2',B_a'$ respectively, $t_0,t_2\in \{0,1\}$ the dimension of the 
locus of such curves is $t_a$, and for the structure of the full automorphism group we have:
\[
G:=\left\{
\begin{array}{ll}
H & \mbox{ if } (n,2)=2, n\equiv 2 \mod 4 \mbox{ and } t_0=t_2=1\\
\mathbb{Z}/n\mathbb{Z}\times S_4 & \mbox{ if in all other cases}
\end{array} 
\right. 
\]
where $H$ is the group given in terms of generators and relations as 
\[
\langle X,Y,T \mid Y^2=X^4=XTX^{-1}=YTY^{-1}=T,
\]
\[
(X^{-1}Y)^3=T^k, T^n=1 \rangle
\]
for some $k \in \{1,\ldots,n\}$.
\end{pro}
\begin{proof}
The dimension argument is clear and the structure follows from \cite[th. 20]{Ko:99}.
\end{proof}

{\bf Case C}
In this case we assume that the extra automorphism group reduces to an element in $S_4$ of order $2$ that is 
conjugate to $\sigma \mu^2$. By computation we find that if we change the coordinate by the M\"obius 
transformation 
$Q(x)={\frac {\sqrt {2}/2\,i \left(  \left( 2+\sqrt {2} \right) x+2-\sqrt {2}
 \right) }{x-1}}
$
the resulting automorphism $Q^{-1}\sigma\mu^2 Q(x)=-x$.
We have the following orbits:
\[
B_0'':=\big\{ Q(P), P\in \{0,\infty,\pm 1,\pm i \}\big\}
\]
\[
B_1'':=\{\mbox{roots of } g_1(Q(x))(x-1)^8, g_1(x):=x^8+14x^4+1 \}
\]
\[
B_2'':=\{ \infty, \mbox{ roots of } g_2(Q(x))(x-1)^{12}, g_2(x):=(x^4+1)(x^8+34x^4+1)\}
\]
\[
B_a'':=\{ \mbox{roots of } f_a(Q(x))(x-1)^{24}, f_a(x):=(x^8+14x^4+1)^3-a (x^5-x)^4\},
\]
Observe, that $B_2''$ can not be in the ramication locus.
\begin{pro}
Assume that  the curve admits an extra automorphism of order $\delta$, so that the image of the extra automorphism 
in the reduced group is conjugate to $\sigma\mu^2$, and that there are $t_0,t_1,t_a$ orbits of the form 
$B_0'',B_1'',B_a''$ respectively, $t_0,t_1\in \{1,0\}$, then the dimension of such curves is equal to $t_a$ and 
for the structure of the full automorphism group we have:
$G=S_4 \rtimes Z_n $
\end{pro}
\begin{proof}
The dimension argument is clear, since we have $t_a$ free variables and the structure follows by \cite[th. 20]{Ko:99}.
\end{proof}
$\mathbf{A_5}$ {\bf Case.}
We will use again the notation of proposition \ref{listcases}.
We observe that every element in $A_5$ is conjugate to $\sigma,\rho,\sigma\rho$ which have orders 
$(5,2,3)$ respectively.
We distinguish the following $3$ cases:

{\bf Case A}. The curve has reduced group $A_5$ and an extra automorphism of 
order $\delta=5$.

The orbit $B_\infty$ can not be included in the ramification 
locus, by the assumption that the curve has an extra 
automorphism of order $\delta$. 
Let $t_0=1$ (resp $t_0^*=1$) if $B_0$ ( resp 
$B_0^\infty$) is included in the ramification locus  and
$t_0=0$ (resp. $t_0^\infty=0$) if not.
Let $t_3$ be the number of orbits of the form $B_a$ that 
are included in the branch locus. The number of 
branch points $s$ equals
\[
s=20t_0+30t_0^\infty + 60 t_3.
\]
\begin{pro}
The locus of curves with reduced group isomorphic to $A_5$ is of dimension 
$t_3$. If $(n,2)=1$ or if $t_0^\infty=0$ then the full automorphism group 
is isomorphic to $A_5 \times \mathbb{Z}/n\mathbb{Z}$. Otherwise the full 
automorphism group admits a complicated presentation in terms of generators
and relations given in \cite{Ko:99}.
\end{pro}
\begin{proof}
The description of the automorphism group is given in theorem \cite[th. 19]{Ko:99}.
For the dimension statement we argue as before: there are $t_3$ free variables.
\end{proof}
{\bf Example:} Let $f_0,f_1$ be as in proposition (\ref{listcases}). 
Using Maple we compute the dihedral invariants of the curve:
\[
y^n=f_1(x)^3-a f_0(x)^5,
\]
\[
\begin{array}{ll}
u_1 = & 2\, \left( -a-684\,i \right) ^{12} \\
u_2 = &-2\,i \left( -a-684\,i \right) ^{10} \left( 55\,a-157434\,i \right) \\
u_3 = &2\, \left( -a-684\,i \right) ^{9} \left( 12527460\,i+1205\,a \right) \\
u_4 = &70\,i \left( -a-684\,i \right) ^{8} \left( 374\,a-2213157\,i \right) \\
u_5 = &2\, \left( -a-684\,i \right) ^{7} \left( -69585\,a-130689144\,i
\right) \\
u_6 = &-2\,i \left( -a-684\,i \right) ^{6} \left( 33211924\,i+134761\,a
 \right) \\
u_7 = &2\, \left( -a-684\,i \right) ^{5} \left( -69585\,a-130689144\,i
 \right) \\
u_8 = &70\,i \left( -a-684\,i \right) ^{4} \left( 374\,a-2213157\,i \right) \\
u_9 = &2\, \left( -a-684\,i \right) ^{3} \left( 12527460\,i+1205\,a \right) \\
u_{10} = &-2\,i \left( -a-684\,i \right) ^{2} \left( 55\,a-157434\,i \right) \\
u_{11} = &2\, \left( -a-684\,i \right) ^{2}\\
u_{12} = &2\\
\end{array}
\]
{\bf Case B}
The curve has reduced automorphism group $A_5$ and an extra automorphism of order $2$.
In this case we compute that the M\"obius transformation 
\[Q(x)={\frac {i \left(  \left( \sqrt {10-2\,\sqrt {5}}-2 \right) x+2+\sqrt {
10-2\,\sqrt {5}} \right) }{ \left( \sqrt {5}-1 \right)  \left( x-1
 \right) }},
\] 
changes $\rho$ so that $Q^{-1} \rho Q(x)=-x$. Moreover the orbits of the action of $A_5$ on this selection of coordinates are 
given by 
\[
B_\infty':=\{Q^{-1}(\infty), \mbox{ roots of } f_0(Q(x))   \left( x-1\right)^{11}\}
\]
\[
B_0':=\{ \mbox{roots of } f_1(Q(x))  \left( x-1\right)^{20} \}
\]
\[
{B_{0}^*}':=\{\infty, \mbox{ roots of } f_2(Q(x))   \left( x-1\right)^{30} \}
\]
\[
B_a':=\{\mbox{ roots of }  \left(  f_1(Q(x))^3 -a f_0(Q(x))^5 \right) \left( x-1\right)^{60}\},
\]
$a\neq -1728i$.
The orbit ${B_{0}^*}'$ contains $\{0,\infty\}$ and can not be in the ramification locus. 
\begin{pro}
Assume that the curve has reduced automorphism group isomorphic to $A_5$ and admits an extra automorphism of order two.
Assume moreover that we have $t_\infty,t_0,t_a$ orbits of the form $B_\infty',B_0',B_a$ respectively.
Then the locus of curves has dimension $t_a$ and the  full automorphism group is $\Z/n\Z \times A_5$.
\end{pro}
\begin{proof}
The dimension of the locus is clear, and the structure follows by \cite[th. 19]{Ko:99}.
\end{proof}
{\bf Case C}
The curve has reduced automorphism group $A_5$ and an extra automorphism of order $3$.
In this case we compute that the M\"obius transformation 
\[
Q_1(x):=\frac{Ax+B}{Cx+D}
\]
where 
\[
A:=-1/2\,{\frac {-\sqrt {-60+36\,\sqrt {5}+30\,i\sqrt {2}\sqrt {5+\sqrt {
5}}-6\,i\sqrt {5}\sqrt {2}E}+3+\sqrt {5}+i\sqrt {2}
E}{\sqrt {-60+36\,\sqrt {5}+30\,i\sqrt {2}\sqrt {5+
\sqrt {5}}-6\,i\sqrt {5}\sqrt {2}E}}},
\]
\[
B:=1/2\,{\frac {3+\sqrt {5}+i\sqrt {2}E+\sqrt {-60+36\,
\sqrt {5}+30\,i\sqrt {2}E-6\,i\sqrt {5}\sqrt {2}
E}}{\sqrt {-60+36\,\sqrt {5}+30\,i\sqrt {2}\sqrt {5+
\sqrt {5}}-6\,i\sqrt {5}\sqrt {2}E}}}
\]
\[
C:=1/12\,{\frac { \left( \sqrt {5}-1 \right)  \left( 5\,\sqrt {3}-\sqrt {
5}\sqrt {3}+i\sqrt {30+6\,\sqrt {5}} \right) }{-5+\sqrt {5}}}
\]
\[
D:=-1/12\,{\frac { \left( \sqrt {5}-1 \right)  \left( 5\,\sqrt {3}-\sqrt 
{5}\sqrt {3}+i\sqrt {30+6\,\sqrt {5}} \right) }{-5+\sqrt {5}}},
\]
and  $E:=\sqrt {5+\sqrt {5}}$., satisfies $Q_1^{-1} \sigma \rho Q(x) =\omega x$,
$\omega$ is a primitive $n$-root of unity.
Moreover the orbits of the action of $A_5$ on this selection of coordinates are 
given by 
\[
B_\infty':=\{Q_1^{-1}(\infty), \mbox{ roots of } f_0(Q(x))   \left( x-1\right)^{11}\}
\]
\[
B_0'':=\{ \mbox{roots of } f_1(Q_1(x))  \left( x-1\right)^{20} \}
\]
\[
{B_{0}^*}'':=\{\infty, \mbox{ roots of } f_2(Q_1(x))   \left( x-1\right)^{30} \}
\]
\[
B_a'':=\{\mbox{ roots of }  \left(  f_1(Q(x))^3 -a f_0(Q_1(x))^5 \right) \left( x-1\right)^{60}\},
\]
$a\neq -1728i$.
The orbit ${B_{0}''}$ contains $\{0,\infty\}$ and can not be in the ramification locus. 
\begin{pro}
Assume that the curve has reduced automorphism group isomorphic to $A_5$ and admits an extra automorphism of order two.
Assume moreover that we have $t_\infty,t_0^*,t_a$ orbits of the form $B_\infty',B_0',B_a$ respectively.
Then the locus of curves has dimension $t_a$ and the  full automorphism group $G$ is given by 
\[
G:=\left\{
\begin{array}{ll}
 \Z/n\Z \times A_5 & \mbox{ if } (n,2)=1 \mbox{ or } t_0^*=0\\
H & \mbox{ otherwise }
\end{array}
\right.
\]
where $H$ is the group defined in terms of generators and relations as follows:
\[
\left\langle X,Y,Z,T | T^n=X^3=1, Y^2=T,Z^2=T (XY)^3=T^\ell, (YZ)^3=T^o \right.
\]
\[
\left. (XZ)^2=T^m,XTX^{-1}=T,ZTZ^{-1}=T, YTY^{-1}=T \right\rangle,
\]
for some integers $m,\ell,o \in \{1,\ldots,n\}$.
\end{pro}
\begin{proof}
The dimension of the locus is clear, and the structure follows by \cite[th. 19]{Ko:99}.
\end{proof}

{\bf The Case} $
\big( \mathbb{Z}/p\mathbb{Z} \times \cdots \times  \mathbb{Z}/p\mathbb{Z}\big) \rtimes \mathbb{Z}/m \mathbb{Z}.
$
We use the notation of  proposition \ref{listcases}. By the assumption on the existence of
an extra automorphism we have that the orbits $B_\infty,B_0$ can not be included in the 
ramification locus of the cover.  Assume that we have $t$ orbits of the form $B_a$. 
We have the following 
\begin{pro}
The locus of curves with reduced automorphism $\mathbb{Z}/p\mathbb{Z} \times \cdots \times \times \big) \rtimes \mathbb{Z}/m \mathbb{Z}$ is of dimension $t$. The automorphism group in this 
case is  
\[
\big(\mathbb{Z}/p\mathbb{Z} \times \cdots \times \mathbb{Z}/p\mathbb{Z} \big) \rtimes \mathbb{Z}/m \mathbb{Z}\big)
\times \mathbb{Z}/n\mathbb{Z}.
\]
\end{pro}
\begin{proof}
The dimension is clearly $t$. For the structure of the whole automorphism 
group we use theorem 12 in \cite{Ko:99}.
\end{proof}

{\bf Projective linear groups}
In this case we assume that the reduced group is isomorphic to a projective linear 
group, $PSL(2,p^t)$ or $PGL(2,p^t)$.
We use the notation of proposition \ref{listcases}. 
We have to distinguish the following cases:

{\bf Case A.}
The extra automorphism of order $\delta$ reduces to a subgroup of the element of order $p^t-1$ ($(p^t-1)/2$ if $G=PSL(2,p^t)$).
Then the orbit $B_\infty$ can not be included 
in the ramification locus. Thus the number of branch points of the cover is given by 
\[
\ell(n+1)(n-1)t_0+t_1 |G|,
\]
where $\ell=1$ (resp. $\ell=\frac{1}{2}$) if $G\cong PGL(2,p^t)$ (resp $G\cong PSL(2,p^t)$) 
and $t_0=1$ (resp.  $t_0=0$)  if the orbit $B_0$ is (resp. not) contained in the branch locus.
\begin{pro}
If there is an extra automorphism of order $\delta \mid p^t-1$ and the reduced 
automorphism group is isomorphic to a projective linear 
group then the locus of such curves is of 
dimension $t_1$. The full automorphism group is isomorphic to 
$\mathbb{Z}/n\mathbb{Z} \times G$ if $(n,2)=1$ or if $t_0$=1. Otherwise the 
structure of the full automorphism group can be described with the 
aid of a restriction map.
\end{pro}
\begin{proof}
The structure part is given in \cite[th. 22]{Ko:99}. For the dimension part we argue as we 
did in the previous cases.
\end{proof}
{\bf Example:} We are using the notation of proposition \ref{listcases}. 
For the curve $p=3$, 
\[
y^n= \left(  \left( {x}^{9}-x \right) ^{8}+1 \right) ^{10}-a \left( {x}^{9
}-x \right) ^{72},
\]
with reduced group $PGL(2,9)$, the dihedral invariants are computed using Maple
\[
u(i)=\left\{
\begin{array}{ll}
2  & \mbox{ if } 9\nmid i,\\
4+4a & \mbox{ if } 9 \mid i.
\end{array}
\right.
\]
{\bf Case B} The extra automorphism $\sigma(x)$ order $\delta$ is reduced to the element of order a divisor of $p^t+1$. 
Then the transformation matrix $Q=-i\frac{x+1}{x-1}$ transforms the coordinate of the projective line so that $\sigma$ is 
an element of the form $\sigma(x)=\zeta x$, $\zeta^\delta=1$.
The orbit $Q^{-1} B_0$ contains now $\{0,\infty\}$ and we can assume that we have $t_\infty$ orbits of the form $Q^{-1} B_\infty$ and 
$t_a$ generic orbits. where $t_\infty=\{\pm 1\}$.

\section{Examples of low genus}
We will describe all curves that are $n$-th covers of the projective line and have genus $\leq 10$ where $n$ is a prime number.
The genus of an $n$-th cyclic cover ramified completely at $s$-points is given by (\ref{eq2}). Since $s\geq 3$ all primes $\geq 19$ give rise 
to curves of genus greater than $10$, we will restrict ourselves to $n=3,5,7,11,13,17,19$. We also observe a necessary condition 
for an automorphism to have an extra automorphism $ n |s $. Therefore the only primes that are of interest to us are $n=3,5$.

$\mathbf{n=3}$. The inequality $g\leq 10$ implies that $s\leq 12$. The multiples of $3$ that give genera $2\leq g \leq 10$
are $s=3,6,9,12$.

$\bullet$ For $s=6$ we have the following curve in normal form:
\[
y^3=x^6+ax^3+1.
\]
The discriminant of the polynomial $x^6+ax^3+1$ is computed to be $3^6(-4+a^2)^3$, so if $a\neq \pm 2$ and the 
characteristic $p\neq 3$, the above curves have genus $4$. The dihedral invariants are computed $u_1=2a^2$, $u_2=2$, 
and the relation $u_{r-1}^{r/2}=2^{r-2} u_1$ holds, so all the above curves admit $D_3 \times \Z/3\Z$ as a subgroup 
of the full automorphism group. Notice that this group equals the normalizer of cyclic Galois cover group in the full automorphism group.

Moreover, when $s=6$ the reduced group $A_4$ in case B can appear, The roots of the polynomial 
\[
f_A:=x^6+\left( 25-15\sqrt{3} \right) x^3 -26+15\sqrt{3},
\]
are exactly the elements of the orbit $B_1$. Therefore the curve
\[
y^3=f_A(x),
\]
has automorphism group 
\[
\langle R,S \mid R^2=S^2,S^6=1, RSR^{-1}=R \rangle.
\]
The normal form of the above curve is computed:
\[
y^3=x^6 +5i\sqrt{2} x^3 +1,
\]
and the dihedral invariants are $u_1=-100,u_2=2$. 

$\bullet$ For $s=9$ we have the following curve in normal form:
\[
y^3=x^9+ax^6+bx^3+1.
\]
The discriminant of the right hand side is computed to be:
\[
D_9=3^9 (27-18 ba - b^2a^2 + 4a^3 +4 b^3)^3,
\]
and if the characteristic $p\neq 3$ and $(a,b)$ is not in the zero locus of $D_9$ then the genus of the above curve is $7$.
Moreover $2n<s$ and the cyclic Galois cover group is normal in the full automorphism group. The dihedral invariants are 
computed: $u_1=b^3+a^3$ and $u_2=2ba$.
If 
\[
u_{r-1}^r -2^{r-2}u_1^2=4b^3a^3-2b^6-2a^6 \neq 0,
\]
{\em i.e.,} if $a/b$ in not a third root of unity, then the full automorphism group is equal to $\Z/3\Z \times \Z/ 3\Z$.
If $a^3=b^3$ then exactly one of the $\pm 1$ are in the branch locus and the automorphism group 
admits the following presentation $\langle R,S \mid R^6=S^3=(RS)^2=1\rangle.$

$\bullet$ For $s=12$ we consider the following curve in normal form:
\[
y^3=x^{12} +a x^9 +b x^6 +cx^3+1.
\]
The discriminant of this curve is computed:
\[
D_{12}=3^{12}(-256-b^2c^2a^2-18c^3ab-18a^3bc+80 b^2ca +192 ca-\]
\[
144a^2b-144c^2b+27a^4+27c^4+128b^2+4b^3a^2+4c^3a^3+4b^3c^2+6c^2a^2-16b^4)^3,
\]
so it the characteristic is not $3$ and if $(a,b,c)$ is not in the zero locus of $D_{12}$ then the 
curve has genus $10$  The dihedral invariants are computed: $u_1=c^4+a^4$, $u_2=c^2b+a^2b$, 
$u_3=2ca$. We compute that 
\[
2^{\frac{r-2}{2} } u_1-u_{r-1}^{\frac{r}{2}}=2(a-c)^2(a+c)^2,
\]
\[
2^{\frac{r-2}{2} } u_1+u_{r-1}^{\frac{r}{2}}=2(a-ic)^2(a+ic)^2,
\]
where $i$ is a primitive $4$-root of unity. Therefore, if $a=\pm c$ then the automorphism group is 
isomorphic to 
\[
\Z /3 \Z \times D_3,
\]
and if $a=\pm ic $ then the automorphism group admits a presentation:
\[
\langle  R,S \mid R^6=S^3=1, RSR^{-1}=S^{-1}\rangle.
\] 
In the case $s=12$ the generic orbit of $A_4$ can appear. The equation of the curve in this case is given by
$y^3=f_a$ where $f_a$ is the polynomial defined in (\ref{eq12}). The corresponding automorphism group equals 
$\Z/2\Z \times A_4$. The dihedral invariants of this curve can be easily computed but the size of the 
expressions prohibits us to express them here.

In the case $s=12$ the $B_2'$ orbit of the group $S_4$ in case $B$ can appear. We compute that the polynomial with roots the 
elements of the orbit $B_2'$ corresponds to the curve:
\[
y^2=g_2(Q(x))(x-1)^8,\;\; g_2(x)=(x^4+1)(x^8+34x^4+1),
\]
where $Q(x)=-1/3 \frac{\big((3+\sqrt{3})x-\sqrt{3}+3 \big) \sqrt{3}}{1-i-(1-i)x}.$ Again, the dihedral invariants are 
computable but too large to be presented here. The automorphism group of this curve is isomorphic to $\Z/3\Z \times S_4$.

\section{An application to fields of definition}

Let $F/K$ be a Galois extension of fields. Let $X$ be a curve defined over $F$, {\em i.e.} there is a  map $X\rightarrow \Spe F$. 
 For every $\sigma \in Gal(F/K)$ we define the curve 
$X^\sigma= X\times_{\Spe F} \Spe F \stackrel{\Spe \sigma}{\longrightarrow} \Spe F$.  

The group $H(X)$ can be defined as 
\[
H(X)=\{\sigma\in Gal(F/K): X^\sigma/\Spe F \cong_F X/ \Spe F\}.
\]
It is known that $H(X)$ is a closed subgroup of 
$Gal(F/K)$ in the Krull topology \cite[prop. 2.1]{Deb-Ems}, and defines the field 
of moduli $F^{H(X)}$. If $[X]$ is the moduli point corresponding
to the curve $X$ then the field of moduli coinsides with the residue field of $M_g[X]$ at 
the point $[X]$ \cite{Baily62}.

A subfield $L$ of $F$ is 
called a field of definition if there is a curve $X'$ defined over $L$ such that
$X,X'$ are isomorphic over $F$, {\em i.e.}, 
$X' \times_K F \cong X$.

It is known that if the  field of moduli is a field of definition 
then it is the smallest field of definition. Otherwise the field of 
definition is a finite extension of the field of moduli. 

We have to notice here that for curves of genera $0,1$ the 
field of moduli and the field of definition coincide. 

Whether or not the field of moduli  is a
field of definition is in general a difficult problem that 
goes back to Weil, Baily, Shimura et al.
Using dihedral invariants we are able to prove the 
following:
\begin{pro}
Consider a cyclic cover $X$ of the projective line given by 
\[
y^n =f(x),
\]
such that $D_\delta$, $\delta \mid n$ is a subgroup of   the reduced 
automorphism group. Then the field of moduli is the field of definition 
and if $u_1\neq 0$, then a  rational model over the field of moduli is given by:
\begin{equation} \label{fieldofmoduli}
y^n=x^{r\delta} + \left(\frac{u_1}{2}\right)^{\frac{1}{r}}x^{(r-1)\delta} + 
\sum_{i=1}^{r-2} \frac{u_{r-i}}{u_1} \left(\frac{u_1}{2} \right)^\frac{r-i}{r} x^{\delta i} +1. 
\end{equation}
\end{pro}
\begin{proof}
Suppose that the curve $X$ has dihedral invariants $\bar{u}:=(u_1,\ldots,u_r)$ and admits a dihedral group 
$D_\delta$ as a subgroup of the full automorphism group. Assume first that $u_1\neq 0$.
We will give a model of the curve $X$ defined over the field of moduli, 
{\em i.e.}, a model of the curve with dihedral invariants $\bar{u}$.
Let us write the curve $X$ in normal form 
\[
y^n=x^{r\delta}+\sum_{i=1}^{r-1} a_i x^{i\delta} + 1. 
\]
Then, by lemma \ref{dihed-sym} we have that $a_{r-i}=a_i$ and by 
proposition \ref{dih-form}
we have that $ 2^{r-2} u_1^2=u_{r-1}$.
Therefore for the coefficients $A_i:=\frac{u_{r-i}}{u_1} \left(\frac{u_1}{2} \right)^\frac{r-i}{r}$
of the right hand side of (\ref{fieldofmoduli}) we have that $A_i$.
The dihedral invariants  $u_i'$ of the curve (\ref{fieldofmoduli}) can now be computed:
\[
u_i'=A_1^{r-i} A_i + A_{r-1}^{r-i}a_{r-i}=2 A_{r-1}^{r-i}A_{r-i}=u_i.
\]
This proves that the dihedral invariants of the curve in (\ref{fieldofmoduli})
are $u_1,\ldots,u_{r-1}$. 

If $u_1=0$ then after a change of dihedral invariants as it explained in Remark 3  on page \pageref{dihedral-blowup}
we can consider suitable dihedral invariants $u_1^{(e)},\ldots,u_r^{(e)}$ so that 
$u_1^{(e)} \neq 0$ and apply the argument we have used for the case $u_1\neq 0$.
\end{proof}

\bibliography{bib.bib}

\end{document}